\documentstyle[11pt]{amsart}
\title{On endomorphisms of surface mapping class groups}
\author{Mustafa Korkmaz}
\date{\today \\ \hspace*{0.3cm}
\\ \hspace*{0.4cm}{\it Key words and phrases:}
Mapping class groups, surfaces, residually finite groups.}
\def\sl{SL_2({\bf Z})}
\begin{document}
\newtheorem{thm}{Theorem}
\newtheorem{cor}[thm]{Corollary}
\newtheorem{prop}[thm]{Proposition}
\newtheorem{lemma}[thm]{Lemma}
\newtheorem{df}[thm]{Definition}
\newtheorem{thmo}{Theorem}
\maketitle
\markboth{MUSTAFA KORKMAZ}{ON ENDOMORPHISMS OF SURFACE MAPPING CLASS GROUPS}
\setcounter{secnumdepth}{2}
\setcounter{section}{0}
\begin{quote}
{\footnotesize ABSTRACT: We prove in this paper that every
endomorphism of the mapping class group of certain orientable
surfaces onto a subgroup of finite index is in fact an
automorphism.}
\end{quote}

\section{Introduction}

Let $S$ be a compact connected orientable surface. The mapping class group
${\cal M}_S$ of the surface $S$ is the group of isotopy classes of
orientation preserving diffeomorphisms $S\to S$. The extended mapping
class group ${\cal M}_S^*$ of $S$ is the group of isotopy classes of all
diffeomorphisms $S\to S$. Note that the isotopy classes of orientation
reversing diffeomorphisms are also included in ${\cal M}_S^*$,
and hence ${\cal M}_S$ is a subgroup of ${\cal M}_S^*$ of index two.

Recall that a group $G$ is called residually finite if for each
$x\neq 1$ in $G$ there exists a homomorphism $f$ from $G$ onto
some finite group such that $f(x)$ is nontrivial. Equivalently,
there is some finite index normal subgroup of $G$ that does not
contain $x$. $G$ is called hopfian if every surjective
endomorphism of $G$ is an automorphism. It is well known that
finitely generated residually finite groups are hopfian \cite{ls}.
$G$ is called cohopfian if every injective endomorphism of $G$ is
an automorphism.

The mapping class group of an orientable surface is finitely
generated \cite{l,b} and residually finite \cite{g,i1}. Hence it
is hopfian. N.V. Ivanov and J.D. McCarthy \cite{im} proved that
${\cal M}_S$ is also cohopfian. Author and J.D. McCarthy \cite{km}
proved that if $\phi :{\cal M}_S\to{\cal M}_S$ is a homomorphism
such that $\phi({\cal M}_S)$ is a normal subgroup and ${\cal
M}_S/\phi({\cal M}_S)$ is abelian, then $\phi$ is an automorphism.

In this paper, we prove further that if $\phi$ is an endomorphism
of the mapping class group ${\cal M}_S$ onto a finite index
subgroup, then $\phi$ is in fact an automorphism, with a few
exceptions. The proof of this result relies on a result of R.
Hirshon, which states that if $\phi$ is an endomorphism of a
finitely generated residually finite group $G$ such that $\phi(G)$
is of finite index in $G$, then $\phi$ restricted to $\phi^n(G)$
is an injection for some $n$.

D.T. Wise \cite{w} gave an example of a  finitely generated
residually finite group $G$ and an endomorphism $\Phi$ of $G$ such
that the restriction of $\Phi$ to $\Phi^n(G)$ is not injective for
any $n$, answering a question of R. Hirshon in negative. It might
be interesting to consider the same question for mapping class
groups of surfaces.

\section{Endomorphisms of mapping class groups}

Let $S$ be a compact connected oriented surface of genus $g$ with
$b$ boundary components. For any simple closed curve $a$ on $S$,
there is a well known diffeomorphism, called a right Dehn twist,
supported in a regular neighborhood of $a$. We denote by $t_a$ the
isotopy class of a right Dehn twist about $a$, also called a Dehn
twist. Note that $ft_af^{-1}=t_{f(a)}$ for any orientation
preserving mapping class $f$.

The pure mapping class group ${\cal PM}_S$ is the subgroup of
${\cal M}_S$ consisting of those orientation preserving mapping
classes which preserve each boundary component.

For a group $G$ and a subgroup $H$ of it, we denote by $C_G(H)$
the centralizer of $H$ in $G$. The center of $G$ is denoted by $C(G)$.

\begin{thm}
Let $G$ be a finitely generated residually finite group, and let
$\phi$ be an endomorphism of $G$ onto a finite index subgroup.
Then there exists an $n$ such that the restriction of $\phi$ to
$\phi^n(G)$ is an injection.
\label{thm1}
\end{thm}

\begin{thm}
Let $S$ be a compact connected orientable surface of genus $g$ with $b$
boundary components. Suppose, in addition, that if $g=0$ then $b\geq 5$, if
$g=1$ then $b\geq 3$, and if $g=2$ then $b\geq 1$. Then any
isomorphism between two finite index subgroups of the extended mapping
class group ${\cal M}_S^*$ is the restriction of an inner automorphism
of ${\cal M}_S^*$.
\label{thm2}
\end{thm}

Theorem \ref{thm1} was proved by R. Hirshon (\cite{h}), and
Theorem \ref{thm2} was proved by N.V. Ivanov \cite{i2} for
surfaces of genus at least two and by the author \cite{k} for the
remaining cases. Since the mapping class group ${\cal M}_S$ is
normal in ${\cal M}_S^*$, we deduce the following theorem.

\begin{thm}
Let $S$ be a compact connected orientable surface of genus $g$ with $b$
boundary components. Suppose, in addition, that if $g=0$ then $b\geq 5$, if
$g=1$ then $b\geq 3$, and if $g=2$ then $b\geq 1$. Then any
isomorphism between two finite index subgroups of the mapping
class group ${\cal M}_S$ is the restriction of an automorphism
of ${\cal M}_S$.
\end{thm}

\begin{lemma}
Let $S$ be a closed orientable surface of genus two and let $\Gamma$ be a
finite index subgroup of ${\cal M}_S$. Then the center $C(\Gamma)$ of
$\Gamma$ is equal to $\Gamma\cap\langle\sigma\rangle$, where $\sigma$ is
the hyperelliptic involution.
\label{lemma}
\end{lemma}
{\it Proof:} Since the subgroup $\langle\sigma\rangle=
\{1,\sigma\}$ is the center
of ${\cal M}_S$, its intersection with $\Gamma$ is contained in the
center of $\Gamma$.

Now let $f\in C(\Gamma)$ and let $N$ be the index of $\Gamma$ in
${\cal M}_S$. Since $t_a^N\in\Gamma$ for all simple closed curves
$a$, we have $t_{f(a)}^N=ft_a^Nf^{-1}=t_a^N$. It follows that
$f(a)=a$ (cf. \cite{im}). Hence, $ft_af^{-1}=t_{f(a)}=t_a$. Since
${\cal M}_S$ is generated by Dehn twists, $f\in C({\cal M}_S)
=\langle\sigma\rangle$. $\Box$

\bigskip
We are now ready to state and prove the main result of
this paper.

\begin{thm}
Let $S$ be a compact connected orientable surface of genus $g$
with $b$ boundary components. Suppose, in addition, that if $g=0$
then $b\neq 2,3,4$, and if $g=1$ then $b\neq 2$. If $\phi$ is an
endomorphism of ${\cal M}_S$ such that $\phi({\cal M}_S)$ is of
finite index in ${\cal M}_S$, then $\phi$ is an automorphism.
\label{mthm}
\end{thm}
{\it Proof:} If $S$ is a (closed) sphere or a disk, then ${\cal
M}_S$ is trivial. Clearly, the conclusion of the theorem holds.

Suppose first that $S$ is a torus with $b\leq 1$ boundary
component. It is well known that ${\cal M}_S$ is isomorphic to
$\sl$. The commutator subgroup of $\sl$ is a nonabelian free group
of rank $2$ and its index in $\sl$ is $12$. Let us denote it by
$F_2$. $\phi(F_2)$ is contained in $F_2$ as a finite index
subgroup. If this index is $k$, $\phi(F_2)$ is a free group of
rank $k+1$. Since there is no homomorphism from $F_2$ onto a free
group of rank $\geq 3$, it follows that $k=1$. That is,
$\phi(F_2)=F_2$. In particular, $\phi(\sl)$ contains $F_2$. The
fact that $\phi$ is an automorphism in this case was proved in
\cite{km}.

Suppose now that $S$ is not one of the surface above and not a
closed a surface of genus $2$. Let us orient $S$ arbitrarily.
Since ${\cal M}_S$ is finitely generated and residually finite,
there exists an $n$ such that the restriction of $\phi$ to
$\phi^n({\cal M}_S)$ is an isomorphism onto $\phi^{n+1}({\cal
M}_S)$. Note that the subgroups $\phi^n({\cal M}_S)$ and
$\phi^{n+1}({\cal M}_S)$ are of finite index in ${\cal M}_S$.
Hence, there is an automorphisms $\alpha$ of ${\cal M}_S$ such
that the restrictions of $\alpha$ and $\phi$ to $\phi^n({\cal
M}_S)$ coincide.

Let $N$ be the index of $\phi^n({\cal M}_S)$ in ${\cal M}_S$.
For any simple closed curve $a$ on $S$, $t_a^N$ is contained in
$\phi^n({\cal M}_S)$. Hence,
$\alpha(t_a^N)=\phi(t_a^N)$. Let $f\in{\cal M}_S$ be any element. Then
\begin{eqnarray*}
t_{f(a)}^N&=&\alpha^{-1}(\phi(t_{f(a)}^N))\\
&=&\alpha^{-1}(\phi(ft_a^Nf^{-1}))\\
&=&\alpha^{-1}(\phi(f))\alpha^{-1}(\phi(t_a^N))\alpha^{-1}(\phi(f^{-1}))\\
&=&\alpha^{-1}(\phi(f))t_a^N\alpha^{-1}(\phi(f))^{-1}\\
&=&t_{\alpha^{-1}(\phi(f))(a)}^N.
\end{eqnarray*}
Hence, $\alpha^{-1}(\phi(f))(a)=f(a)$ for all $a$ (cf. \cite{im}).
It follows that $f^{-1}\alpha^{-1}(\phi(f))$ commutes with all
Dehn twists. Since ${\cal PM}_S$ is generated by Dehn twists, it
is in $C_{{\cal M}_S}({\cal PM}_S)$, the centralizer of ${\cal
PM}_S$ in ${\cal M}_S$. But $C_{{\cal M}_S}({\cal PM}_S)$ is
trivial \cite{im}. Hence, $\alpha^{-1}(\phi(f))=f$.  Therefore,
$\phi=\alpha$. In particular, $\phi$ is an automorphism.

Suppose finally that $S$ is a closed surface of genus two. Let $R$
be a sphere with six holes. Then ${\cal M}_R$ is isomorphic to the
quotient of ${\cal M}_S$ with its center $\langle\sigma\rangle$,
where $\sigma$ is the hyperelliptic involution (cf. \cite{bh})  .
Let us identify ${\cal M}_R$ and ${\cal M}_s/ \langle
\sigma\rangle$, and let $\pi: {\cal M}_S\to {\cal M}_R$ be the
quotient map. Since $\phi(\sigma)$ is in the center of $\phi({\cal
M}_S)$, either $\phi(\sigma)=\sigma$ or $\phi(\sigma)=1$ by the
lemma above.

If $\phi(\sigma)=\sigma$, then $\phi$ induces an endomorphism
$\Phi$ of ${\cal M}_R$, such that $\pi\phi=\Phi\pi$. Then we have
a diagram in which all squares are commutative:
\[
\begin{array}{ccccccccc}
1& \longrightarrow &\langle\sigma\rangle& \longrightarrow&{\cal M}_S&
\stackrel{\pi}\longrightarrow&{\cal M}_R&\longrightarrow &1\\
&  &\Big\downarrow {\rm I}& &\Big\downarrow \phi& &\Big\downarrow\Phi& & \\
1& \longrightarrow &\langle\sigma\rangle& \longrightarrow&{\cal M}_S&
\stackrel{\pi}\longrightarrow &{\cal M}_R & \longrightarrow &1
\end{array}
\]
where ${\rm I}$ is the identity homomorphism. Since the image
$\Phi({\cal M}_R)$ of $\Phi$ is of finite index, $\Phi$ is an
automorphism by the first part. By $5$-lemma, $\phi$ is an
automorphism.

If $\phi(\sigma)=1$, then $\phi$ induces a homomorphism
$\overline{\phi}:{\cal M}_R\longrightarrow {\cal M}_S$ such that
$\overline{\phi}\pi=\phi$.
 The image of the endomorphism
$\Phi=\pi\overline{\phi}$ of ${\cal M}_R$ has finite index. Since
$R$ is a sphere with six holes, $\Phi$ is an automorphism by the
first part. Then, $\phi$ is an automorphism, and hence $\sigma=1$.
This contradiction finishes the proof of our theorem. $\Box$

\bigskip
\noindent {\bf Remark}: If $S$ is a sphere with two holes, then
${\cal M}_S$ is a group of order two, and if $S$ is a sphere with
three holes, then ${\cal M}_S$ is isomorphic to the symmetric
group on three letters. Hence, in these cases the trivial
homomorphism is an endomorphism onto a finite index subgroup which
is not automorphism. We do not know if the conclusion of Theorem
\ref{mthm} holds if a sphere with four holes and a torus with two
holes.

\bigskip
{\sc Department of Mathematics, Middle East Technical University,
06531 Ankara, Turkey. E-mail: {\tt
korkmaz\verb+@+math.metu.edu.tr.}}

\end{document}